\documentclass{amsart}
\usepackage{pdfpages,  
latexsym,amssymb}
\def\proof{\noindent{\sc Proof. }}
\def\EOP{\hfill$\Box$}
\def\natnum{\hbox{\rm I\kern-.17em N}}
\def\integ{\hbox{\rm Z\kern-.3em Z}}
\def\disj{\uplus }

\newtheorem{thm}{Theorem}[section]

\newtheorem{prop}{Proposition}[section]
\newtheorem{defn}{Definition}[section]
\newtheorem{cor}{Corollary}[section]
\newtheorem{qn}{Question}[section]
\newtheorem{examp}{Example}[section]

\newtheorem{rk}{Remark}[section]

\begin{document}
\title{Discrete Morse functions from lexicographic orders}
\author{Eric Babson}
\address{Department of Mathematics, Box 354350\\
        University of Washington\\
        Seattle, WA 98195}
\email{babson@math.washington.edu}
\author{Patricia Hersh}
\address{Department of Mathematics\\
University of Michigan\\
2074 East Hall\\
525 East University Ave.\\
Ann Arbor, MI 48109-1109}
\email{plhersh@umich.edu}
\maketitle

\begin{abstract}
This paper shows how to construct a discrete Morse function with a relatively
small number of critical cells for the order complex of any finite poset
with $\hat{0} $ and $\hat{1}$ from any lexicographic order on its maximal
chains.  Specifically, if we attach facets according to the lexicographic
order on maximal chains, then each facet contributes at most one new
face which is critical, and at most one Betti number changes;
facets which do not change the homotopy type also
do not contribute any critical faces.  Dimensions of
critical faces as well as a description of which facet attachments change
the homotopy type are provided in terms of interval systems 
associated to the facets.  As one application, the 
M\"obius function may be computed as 
the alternating sum of Morse numbers.

The above construction enables us to prove
that the poset $\Pi_n/S_{\lambda }$ of partitions of 
a set $\{ 1^{\lambda_1 },\dots ,k^{\lambda_k }\} $ with repetition is 
homotopy equivalent to a wedge of spheres of top dimension
when $\lambda $ is a hook-shaped partition; it is likely that the proof
may be extended to a larger class of $\lambda $ and 
perhaps to all $\lambda $, despite a result of Ziegler which shows that 
$\Pi_n/S_{\lambda }$ is not always Cohen-Macaulay [Zi].  Additional
applications appear in [He2] and [HW].
\end{abstract}

\section{Introduction}
\begin{figure}[h]
\begin{picture}(250,100)(-75,10)
\includegraphics[width=0.37\textwidth, angle = 0]{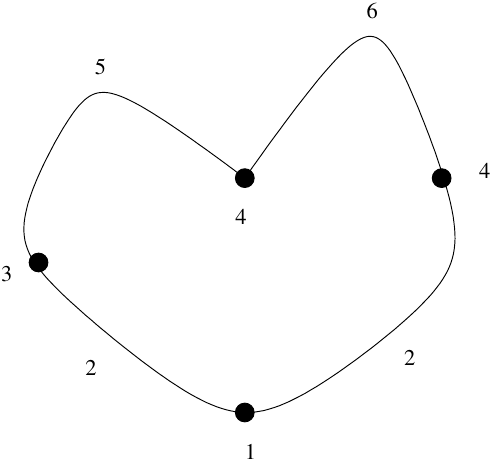}
\end{picture}
\caption{A discrete Morse function}
\label{height}
\end{figure}

Let us begin by reviewing
the notion of a discrete Morse function, as established by 
Robin Forman in [Fo].  A function $f$ which 
assigns real values to the cells in a 
CW-complex is a discrete Morse function if (1) for each $d$-cell 
$\sigma $, denoted $\sigma^{(d)}$, the
sets $ \{ \tau^{(d-1)} \subseteq \overline{\sigma } | f(\tau^{(d-1)} ) 
\ge f(\sigma )\} $ and $ \{ \overline{\tau^{(d+1)} }\supseteq \sigma
| f(\tau^{(d+1)} ) \le
f(\sigma )\} $ each have cardinality at most one, and (2) $f(\sigma )
\ge f(\tau )$ for $\sigma^{(d)} \subseteq \overline{\tau^{(d+1)} }$ 
implies that $\sigma $ is a regular face 
of $\tau $.  Requiring condition (1) for every
$\sigma $ implies that 
for each $\sigma $, at most one of these cardinalities is greater than 0;
when both cardinalities are 0, then $\sigma $ is called 
a {\bf critical cell}.
Figure ~\ref{height} 
gives an example of a discrete Morse function on a 1-sphere resulting from
a height function.

Critical cells record changes in topological structure as a complex is 
built by sequentially inserting cells in the order specified by the Morse
function, with critical $i$-cells behaving similarly to critical points 
of index $i$ in a traditional Morse function; the 
noncritical cells may be eliminated by elementary 
collapses without changing the
homotopy type.  This collapsing results from the fact that
noncritical cells come in pairs 
$\sigma^{(p)} \subseteq \overline{\tau^{(p+1)}} $ 
which prevent each other from being 
critical by satisfying $f(\sigma )\ge f(\tau )$.  Thus, a 
discrete Morse function $f$ gives rise to 
a matching on the noncritical cells in the {\bf face poset}, namely in the 
partial order on cells defined by letting $\sigma $ be 
less than $\tau $ for each $\sigma $ in the boundary of $\tau $ .
\begin{defn}
A matching is {\bf acyclic} if the directed 
graph obtained by directing face poset edges belonging to the matching
upwards and directing all other poset edges downward has no directed 
cycles.
\end{defn}  

Observe that the face poset 
matching resulting from a discrete Morse function is always
acyclic, since the edges are oriented
in the direction in which $f$ decreases.  Conversely,
many different (but in some sense equivalent) discrete Morse functions 
may be constructed from any face poset acyclic
matching.  The unmatched cells in an acyclic matching
are exactly the critical cells of any corresponding
discrete Morse function.  
We shall work primarily in terms of acyclic matchings, but with
this correspondence (due to Chari [Ch]) always in the background.

Let us deviate slightly from Forman's conventions by including the 
empty set in the domain of $f$ and in the face posets upon which we construct
acyclic matchings; this has the effect of expressing results in terms of
reduced Euler characteristic and reduced homology.
Let $m_i$ be the number of critical cells of dimension $i$ in a discrete
Morse function, let $b_i$ be the Betti number recording the rank of 
$H_i(\Delta )$, and let 
$M_a$ be the partial complex comprised of those cells $\sigma $ such that 
$f(\sigma )\le a$.  When using the reduced version of discrete Morse 
theory, we denote the number of critical $i$-cells by $\tilde{m}_i$ and
the reduced Betti numbers by $\tilde{b}_i$.
It is shown in [Fo] that $\Delta $ collapses onto a complex 
$\Delta^M$ comprised of the discrete Morse function critical cells in a
way that preserves homotopy type.  Specifically, Forman shows that 
the following results carry over from traditional Morse theory, the
first two of which are called the Morse inequalities.
\begin{enumerate}
\item
$\tilde{m}_j \ge \tilde{b}_j $ for $-1 \le j\le dim(\Delta )$
\item
$\sum_{i=0}^{j+1}
(-1)^i \tilde{m}_{j-i} 
  \ge \sum_{i=0}^{j+1} (-1)^i \tilde{b}_{j-i}$ 
  for $0\le j\le dim(\Delta )$, with 
equality achieved when $j=dim(\Delta )$
\item
If there are no critical cells $\sigma $ satisfying $a \le f(\sigma )\le b$ 
for some $a<b$, then the partial complex $M_b$ collapses onto the 
partial complex $M_a$.  In particular, $\Delta $ is homotopy 
equivalent to the complex resulting from a series
of such collapses.
\end{enumerate}
From these facts, it immediately follows that:
\begin{enumerate}
\item
If $\tilde{m}_i=0$ for all $i$ then $\Delta $ is collapsible.
\item
If $\tilde{m}_i=0$ for all $i\ne j$ for some fixed $j$, then $\Delta $ is 
homotopy equivalent to a wedge of $j$-spheres.
\end{enumerate}

Figure ~\ref{height}
gives an example of a (non-reduced) discrete Morse function with
$b_0=b_1 =1$ and $m_0=m_1=2$.  Letting $f(\emptyset )= 1.5$ turns this into
a reduced discrete Morse function with $\tilde{b}_0=0, \tilde{m}_0=
\tilde{b}_1=1$ and $\tilde{m}_1=2$.  The Morse 
numbers are larger than the Betti numbers because there is a critical cell
of dimension 0 that is labeled 4 which locally looks as 
though it is creating a new connected
component as the complex is built from bottom to top and there is a 
critical cell of dimension one that is labeled 5 which locally appears to be 
closing off a 1-cycle, but these two critical cells actually cancel each
other's effect.  Forman notes that whenever an acyclic
matching has two critical cells $\sigma^{(p)} $ and $\tau^{(p+1)}$ such 
that there is a unique gradient path from $\tau $ to $\sigma $ (i.e. a path
upon which $f$ decreases at each step), then one
may obtain a new acyclic matching in which $\sigma $ and $\tau $ are no 
longer 
critical by reversing this gradient path, since doing so cannot create any
directed cycles.  Such 
a reversal for instance straightens our 1-sphere into a standard circle, 
which
only has one critical 0-cell and a single critical 1-cell.
This mechanism for cancelling critical cells
will be particularly useful in conjunction with
``lexicographic discrete Morse functions'' because they are well-suited to 
verifying uniqueness of gradient paths, as demonstrated in the application
of Section 5.

We shall consider finite posets with unique minimal and maximal
elements $\hat{0} $ and $\hat{1}$.  A chain $\hat{0} = u_0 < u_1 < \cdots
< u_r = \hat{1}$ of comparable poset elements is 
{\bf saturated} (or {\bf maximal}) if 
$u_i \le v \le u_{i+1}$ implies $v=u_i$ or $v=u_{i+1}$ for $0\le i <r$.  
A poset $P$
is {\bf graded} if for each $x\in P$, all saturated chains from 
$\hat{0}$ to $x$
have the same length, called the {\bf rank} of $x$.  We will not require
our posets to be graded, but we will refer 
to the rank of an element within a 
saturated chain, not requiring this to 
agree with its rank within other saturated chains.

\begin{figure}[h]
\begin{picture}(250,110)(-30,10)
\includegraphics[width=0.6\textwidth, angle = 0]{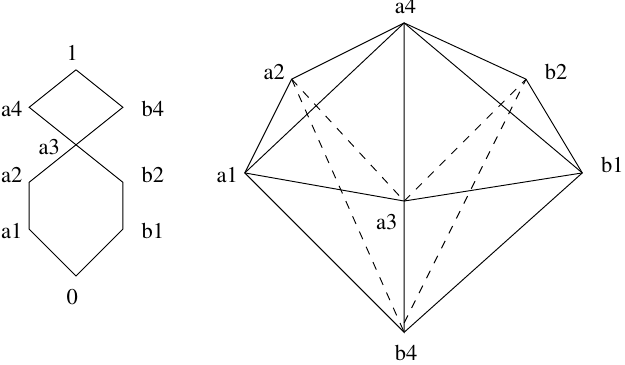}
\end{picture}
\caption{A poset and its order complex}
\label{delta_p}
\end{figure}

Recall that the {\bf order complex}, denoted $\Delta(P)$, for a 
finite poset $P$ with $\hat{0}$ and $\hat{1}$
is the simplicial complex comprised of an 
$i$-face for each chain $\hat{0} < v_0 < \cdots < v_i <\hat{1}$ of 
$i+1$ comparable poset elements.  For instance, the order complex of the 
face poset of any simplicial complex $K$ is the first
barycentric subdivision of $K$, allowing one to apply poset
methods to arbitrary simplicial complexes (for example, Kozlov uses 
this approach on certain graph complexes in [Ko2]).  Figure ~\ref{delta_p}
provides an example of a poset and its order complex.  
Our discrete Morse function will use an ordering on facets of $\Delta (P)$, 
i.e. on saturated chains in $P$.

The {\bf M\"obius function} of a finite poset $P$ is a function
on pairs of elements, defined as follows:
$\mu_P(u,u)=1$ and $\mu_P(u,w) = -\sum_{u\le v < w} \mu_P(u,v)$; M\"obius
functions arise in combinatorics as the coefficients in
inclusion-exclusion counting formulas.  The 
M\"obius function of $P$ satisfies 
$\mu_P (\hat{0},\hat{1}) = \tilde{\chi }(\Delta (P) ) $ (cf. [Ro]), so
our discrete Morse functions for poset order complexes will give a way
to compute $\mu_P(x,y)$ as the alternating sum of Morse numbers on
the order complex restriction $\Delta (x,y)$, as discussed in Section 3.

\begin{rk}
We will sometimes say that a poset $P$
possesses a topological or geometric 
property (such as shellability), by which we mean
that its order complex possesses this property.
\end{rk}

Section 2 gives a way of constructing a discrete Morse function on the 
order complex of any finite poset with $\hat{0} $ and $\hat{1} $ from any 
lexicographic order on its saturated chains.  In particular, we 
provide an acyclic matching on 
the face poset of the order complex of $P$ (not to be confused with the 
original poset $P$).  To this end, let us partition the set
of faces in $\Delta(P)$ by assigning each face
to the earliest facet containing it, i.e.
by the partitioning
$$ \Delta(P) = F_1 \disj (F_2 \setminus F_1) \disj \cdots \disj 
(F_k \setminus 
(\cup_{j<k}F_j ))$$
where $F_1,\dots ,F_k$ is a lexicographic order (as defined below)
on the facets of $\Delta(P)$ or equivalently on the saturated chains of $P$.
We call the above decomposition based on a lexicographic order a 
{\bf lexicographic decomposition} of the order complex.  Our approach
will be to give a matching on the set of faces in each piece $F_j \setminus 
(\cup_{i<j} F_i)$ of the decomposition such 
that (1) the union of these matchings is an acyclic matching for 
all of $\Delta(P)$ and (2) each $F_j \setminus (\cup_{i<j} F_i )$ 
contributes at most one critical cell. 
We will utilize structure on $F_j \setminus (\cup_{i<j} F_i )$
resulting from the use of a 
lexicographic order on facets
in much the way that lexicographic shellability works, but without requiring
our posets to be shellable.  (See Section ~\ref{shell_section} 
for the connection between
discrete Morse functions and shellings.)  In [BH], 
Billera and Hetyei made use of
this same structure resulting from a poset lexicographic order 
in their characterization of the bounding hyperplanes for the 
cone of flag vectors for the set of all graded posets.

\begin{defn}
A poset {\bf lexicographic order} is a total ordering of facets with the 
following property: let two facets $F_1,F_2$ share a face $\sigma $
of ranks $1,2,\dots ,i$, let $\tau\subseteq F_1$ and $\mu\subseteq
F_2$ be faces of ranks $1,2,\dots ,i+1$ such that $\tau\ne \mu $.  If
$F_1$ precedes $F_2$, then any facet containing 
$\tau $ must come before any facet containing $\mu $.  
\end{defn}

Lexicographic orders may be obtained by labeling Hasse diagram 
edges with positive integers and then lexicographically ordering the 
sequences of labels which are associated to the saturated chains.  Such a
labeling is called an {\bf edge-labeling}, and when this gives rise to a 
certain type of shelling (see Section 4) which was introduced 
in [Bj1], the labeling is called an {\bf EL-labeling}.
Recall that the Hasse diagram of $P$ is the graph whose vertices are 
elements of $P$ and whose edges are pairs $(u,w)$ of comparable poset 
elements with no intermediate comparable elements, i.e. pairs satisfying 
$u<w$ such that there is no intermediate $v$ 
satisfying $u<v<w$.  We denote such a minimal comparability by 
$u\prec v$ and call this a {\bf covering relation}.  Many important
classes of posets possess EL-labelings (e.g. see
[St1] for supersolvable lattices, [Bj1] for geometric lattices (a class
of posets which includes all intersection lattices of central hyperplane 
arrangements, and so gives information about the cohomology of the 
complement of such an arrangement), or 
see Bj\"orner's survey article [Bj3, Section 11] for 
additional references).

It is easy to show that any lexicographic order may be achieved by a 
slight generalization of edge-labeling, known as chain-labeling.  
A {\bf chain-labeling} is a labeling of Hasse diagram edges in which the 
label assigned to an edge $u\prec v$ is allowed to depend on the choice of 
rooted chain
$\hat{0}\prec u_1\prec\cdots u_k = u$ leading up to the covering relation
$u\prec v$.  Chain-labeling was introduced in [BW] as part of a tool 
called {\bf CL-shellability}, which generalized EL-shellability and
enabled Bj\"orner and Wachs
to prove that the Bruhat order of any Coxeter group is shellable (and in
particular that each interval is homeomorphic to a sphere).  For other
interesting examples of CL-shellable posets, see for instance [We] for the 
poset of direct sum decompositions of a finite-dimensional
vector space over $F_q$, or [Sh] for the 
subgroup lattice of a finite, solvable group.  From our viewpoint, 
EL-labelings and CL-labelings will give rise to particularly well-behaved
lexicographic discrete Morse functions: it is well-known that such a labeling
for a poset $P$ implies that each interval in $P$ 
is collapsible or homotopy equivalent to a wedge of spheres of 
top dimension, with these
spheres indexed by the saturated chains whose label sequences are 
strictly decreasing in value at each step.

\begin{rk}
For convenience, we will 
also require that any two covering
relations $u\prec v$ and $u\prec w$ with distinct elements $v\ne w$ covering 
the same element $u$ (in the case of an edge-labeling) or covering the 
same root $\hat{0} \prec u_0 \prec \cdots \prec u_k = u$ (in the case of a 
chain-labeling) be assigned distinct labels, 
so that the edge-labeling or chain-labeling gives a  
total order on saturated chains.
\end{rk}

Translating from acyclic matchings to real-valued functions, a lexicographic 
discrete Morse function will assign smaller values to faces belonging to 
lexicographically earlier facets.  Each time a facet is inserted in 
lexicographic order, at most
one of the new faces it adds to the complex is critical, so 
the Betti numbers cannot change very much each
time a facet is inserted.  If the homotopy type changes, 
then either a new cycle is created at
some dimension, or an old cycle becomes a boundary.  
This reflects the fact that each facet
overlaps with the union of earlier ones in either a sphere or a disk, by
virtue of our use of a lexicographic order on facets (cf. [BH] or 
Proposition ~\ref{overlap}).

In Section 5, we use a chain-labeling to construct a 
lexicographic order on the poset $\Pi_n/S_{\lambda }$ of 
partitions of a multiset.  Critical cells in the resulting lexicographic
discrete Morse function live in various dimensions, and the sequence of
critical cell dimensions, listed in the order that they are
inserted in the complex, is not a 
non-increasing sequence.  Thus, we cannot
conclude directly from our original Morse function that 
$\Pi_n/S_{\lambda }$ has the homotopy type of a wedge of spheres.  
However, for certain classes of $\lambda $ we are able to cancel all
the lower-dimensional critical cells in order to show that the order 
complex is homotopy equivalent to a wedge of spheres of 
top dimension.  We cancel critical cells in
pairs $\tau^{(p+1)},\sigma^{(p)}$ by proving there is a unique gradient path
from $\tau $ to $\sigma $.  
The critical cells which we pair together are chosen to 
agree up until a high rank,
forcing any gradient path from $\tau $ to $\sigma $ to preserve
all the lower ranks where the cells agree, aiding us in proving uniqueness.
Our main tool for proving uniqueness is then 
the use of a permutation statistic
based on inversions in the label sequences assigned to the saturated chains;
each step in a gradient path eliminates a single inversion allowing us to 
associate reduced expressions for permutations to gradient paths and then to
use the fact that a cyclic shift permutation has a unique reduced expression
and that its indices are decreasing
to conclude that an associated gradient path is unique.  We hope that our
approach to this example may also be useful for other posets, particularly 
when they have the
homotopy type of a wedge of spheres concentrated in a single dimension, 
whether or not they are Cohen-Macaulay.

Ziegler previously showed that the poset of partitions of an integer (a 
special case of our multiset partition poset) is not Cohen-Macaulay and 
hence not shellable [Zi]; lexicographic discrete Morse functions relax
the requirements of shellability, enabling us to prove that $\Delta (\Pi_n
/S_{\lambda })$ is homotopy equivalent to a wedge of spheres of
top dimension when $\lambda $ is of hook-shape (see 
Section ~\ref{homotopy} for definitions).  

\begin{qn}
Is the poset of partitions of a multiset always homotopy equivalent to
a wedge of spheres of top dimension?  Does
our proof for $\lambda $ of hook-shape extend to all $\lambda $, i.e. is
it true that the puzzle in Section 5.8 has no solutions?
\end{qn}

Examples of lexicographic discrete Morse functions for several other
posets appear in [He2] and [HW].  All of these examples use very natural 
edge-labellings, and they lead us to believe that the most natural of
labellings will often yield Morse functions which may be transformed by
gradient path reversal into perfect Morse functions, i.e. ones in which
the Morse numbers equal the Betti numbers.  Thus, the original complex
may be collapsed onto a cell complex which computes its homotopy type,
in that each cell actually contributes to one of the Betti numbers.

\section{Lexicographic discrete Morse functions}

Given a lexicographic order $F_1,\dots ,F_k$
on the facets in the order complex of a graded poset with $\hat{0} $ and 
$\hat{1}$, we specify a resulting matching on the faces in
$F_j \setminus (\cup_{i<j} F_i)$ for $1\le j\le k$ such that the union of 
these matchings will be an acyclic matching on $\Delta(P)$.  Each facet 
contributes at most one critical cell; it contributes one when the facet
insertion changes the homotopy type of the complex
and does not contribute any when
$\cup_{i\le j} F_i$ may be collapsed onto $\cup_{i<j} F_i $.

This matching will 
require some notation.  If $F_i$ precedes $F_j$ in the lexicographic order
on facets, denote this by $F_i <_{lex} F_j$.  Denote the lexicographically 
earliest facet containing a face $g$ by $eg$, since $eg$ may be thought of 
as an extension of a poset chain $g$ to a saturated chain $eg$.  
If $eg <_{lex} eh$ then we also say $g <_{lex} h$, and it will turn out that
$f(g) < f(h)$ for every discrete Morse function $f$ which is consistent with 
the acyclic matching we are about to construct.
If a face $a$ is matched with a face $b$, denote this by $a\sim b$.  
We will specify the matching
by describing it for the fibres $e^{-1} (F_j)$ for each saturated 
chain $F_j$.

\begin{rk}
Given a saturated chain $F_j$ and a lexicographic order $F_1,\dots ,F_k$ on 
the facets, each maximal face in $F_j\cap (\cup_{i<j} F_i ) $ consists of 
a subchain of $F_j$ given by skipping a single interval of consecutive
ranks in $F_j$.
\end{rk}
  
We refer to these intervals in ranks as 
the {\bf minimal skipped intervals} of $F_j$, and 
we denote this set of intervals by $I(F_j)$.  We call the number of ranks
which are skipped to obtain a maximal face in 
$F_j\cap (\cup_{i<j} F_i ) $ the {\bf height} of the interval.
A minimal skipped interval is called {\bf nontrivial} if it skips more than
one rank.  

\begin{rk}
One way to determine the minimal skipped intervals for a particular $F_j$
is by considering each covering relation $u\prec v$ in $F_j$ in turn; if 
there is
some $u\prec v'$ which comes lexicographically earlier, then one obtains
a face in $F_j\cap (\cup_{i<j} F_i ) $ by taking the intersection of 
$F_j$ with any facet $F_i$ which includes $u\prec v'$, agrees with $F_j$
below $u$, and agrees with $F_j$ above $w\in F_j$ for some $w>v'$ of 
smallest possible rank in $F_j$.  All
maximal faces in $F_j\cap (\cup_{i<j} F_i ) $ arise this way.
Natural labellings often seem to allow simple
classification of minimal skipped intervals. 
\end{rk}

In [BH], the collection $I(F_j)$ of minimal skipped intervals for $F_j$ is
called the {\bf interval system} of $F_j$, and is used to determine the 
bounding hyperplanes for the cone of all possible flag $f$-vectors for
graded posets.  We shall also use interval systems, in our case to
describe the faces in $F_j\setminus (\cup_{i<j}F_i)$, so as to construct
an acyclic matching on the set of faces in each 
$F_j\setminus (\cup_{i<j}F_i)$.  Thus, our
lexicographic discrete Morse function construction
will also work for balanced
simplicial complexes which have the above property of 
poset lexicographic orders that
every maximal face in $F_j\cap (\cup_{i<j} F_i)$ skips a single interval.
See [He1] for a notion of lexicographic order on balanced simplicial
complexes.  
Not all lexicographic orders on simplicial complexes have such an 
interval system structure.

\begin{figure}[h]
\begin{picture}(260,140)(0,10)
\includegraphics[width=0.7\textwidth, angle = 0]{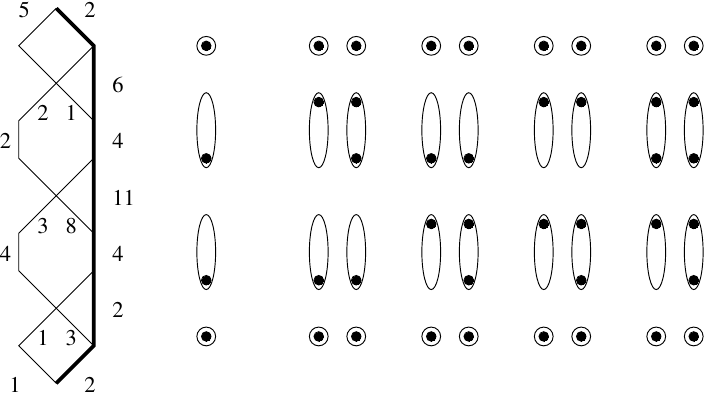}
\end{picture}
\caption{A matching example}
\label{matching}
\end{figure}

Note that a face in
$F_j$ also belongs to a lexicographically earlier facet if it entirely 
misses any of the minimal skipped intervals of $F_j$, so the faces in 
$F_j\setminus (\cup_{i<j} F_i)$ are those faces within $F_j$
which hit every minimal 
skipped interval of $F_j$.  Our matching will be specified in terms of
how a closely related set of nonoverlapping intervals called 
$J(F_j)$ are hit.  Figure ~\ref{matching}
provides an example in which $I(F_j)=J(F_j)$.
The skipped intervals in $I(F_j)$ for the 
darkened saturated chain $F_j$ have 
rank sets $\{ 1\} , \{ 2,3\} , \{ 4,5 \} , \{6 \} $
and the faces in $F_j \setminus (\cup_{i<j} F_i )$ 
are depicted to the right of the poset by showing for each face which 
ranks it hits in each of the intervals in $I(F_j)$.  
The critical cell is listed first, followed by pairs of critical cells
which are to be matched.

Now let us informally describe the matching.
If some vertex $v\in F_j$ does not belong to any 
minimal skipped interval of $F_j$, then including or excluding $v$ has
no bearing on whether a face of $F_j$ also belongs to a lexicographically
earlier facet, so we choose the lowest rank such $v\in F_j$ and match
faces including $v$ with those excluding it.  In this case, the matching on
$F_j \cap (\cup_{i<j} F_i )$ has no critical cells.
Otherwise, every vertex in $F_j$ belongs to some 
minimal skipped interval, and 
then we truncate the intervals in $I(F_j)$
to obtain a collection $J(F_j)$ of nonoverlapping intervals
as follows: we include in $J(F_j)$ the lowest rank interval from $I(F_j)$
in its entirety, truncate all the remaining intervals
to avoid overlap with this interval, discard any intervals which are no
longer minimal and then iterate the procedure on the truncated intervals
which are still minimal among the set of skipped intervals (i.e. 
which do not contain any of the other intervals) to obtain a collection $J$
of nonoverlapping intervals.  

If the $J$-intervals still cover all ranks, then
the unique critical cell in the $F_j$-piece of the lexicographic decomposition
will consist of the minimal
elements of all the intervals in $J(F_j)$.  If the $J$-intervals leave
$i$ ranks uncovered, then there are $2^i$ cells which include exactly
the minimal element for each $J$-interval.  
The other cells in $F_j \setminus (\cup_{i<j} F_i )$
are matched based on the first interval in $J(F_j)$ 
where they deviate from these $2^i$ cells.  This deviation comes either 
by leaving out the minimal
element of the interval or else by including other vertices along with the 
minimal element; we match a face which omits the 
minimal element in the interval with a face which has just that vertex added.
Notice that the faces in $F_j \setminus (\cup_{i<j} F_i )$ need not hit all
the $J$-intervals, but our matching relies on the fact that each such face
does necessarily include an element from the first $J$-interval where it
deviates from the potential critical cell.  To see this, notice that
this first 
$J$-interval is contained in some $I$-interval which may extend downward 
below the $J$-interval, and that the face must hit this $I$-interval; the
face cannot hit the portion of the $I$-interval which is truncated to obtain
the $J$-interval, because we assume that the face agrees with the critical 
cell below the $J$-interval in question.  
The $2^i$ remaining unmatched cells are matched by including/excluding 
the lowest $v\in F_j$ which is not covered by any $J$-interval. 

It may happen that the 
intervals in $I(F)$ cover the saturated chain even though the intervals
in $J(F)$ do not cover it, in which case our construction does not produce
a critical cell.  As an example, consider a poset of rank 5 and a facet
with minimal skipped
intervals $[1,2],[2,3],[3,4]$ in $I(F)$, which in turn yields truncated 
intervals $[1,2],[3,3]$ for $J(F)$ since $[3,4]$ strictly contains the interval
$[3,3]$ and thus is discarded.  Note that the intervals in $J(F)$ do not cover
the vertex of rank 4, enabling $\{ 1,3 \}$ to be matched with $\{
1,3,4\} $, so that no critical cell remains.

Let us formalize the matching construction for any particular $F_k$.
If the $I$-intervals do not cover all ranks, let $j_0$ be the set of
uncovered elements, and let $\rho_0$ denote the lowest rank element
whose
rank is in $j_0$; in this case, match $h\sim h\setminus \{ \rho_0 \} $
for each face $h\in e^{-1}(F_k)$ containing $\rho_0$.  Otherwise,
we derive from the set $I(F_k)$ 
of minimal skipped intervals a collection $J(F_k)$ of disjoint 
intervals by the following algorithm (which is initialized to
$J=\emptyset $):
\begin{enumerate}
\item
Add to $J$ the interval $(u,v)$ in 
$I$ with $u$ of smallest possible rank.
\item
Replace $I$ by the restriction of $I$ to ranks above $rk(v)$.
\item
Delete from $I$ those skipped intervals which are no longer minimal.
\item 
Repeat these steps until $I=\emptyset $.
\end{enumerate}
This algorithm outputs a collection $J = \{ j_1,\dots ,j_r\} $ 
of disjoint intervals.


Let $\rho_i$ denote the lowest rank element from interval
$j_i\in J$ for $1\le i \le r$. 
Let $j_{r+1 }$
be the set of elements not in any $J$-interval, and let $\rho_{r+1}$
denote the lowest rank element in $j_{r+1}$.  Now define a map
$$\tau : e^{-1} (F_k) \rightarrow [r]\cup \{ \infty \} $$ from chains $h$
that extend to $F_k$ to indices of $J$-intervals by
$$h\mapsto \min_{1\le i\le r} \{ h\cap j_i \ne \{ \rho_i \} \} ,$$
letting $\tau h = \infty $ when the set $\{ i | 1\le i \le r,
h\cap j_i \ne \{ \rho_i \} \} $ is empty.
The matching is defined by $h\sim h \backslash \{ \rho_{\tau h} \} $ for
$\rho_{\tau h}\in h$ and $h\sim h \cup \{\rho_{\tau h} \}$ otherwise.
This leaves $h$ not yet matched when $\tau h = \infty $, namely when $j_0$ is
empty and $h$ is comprised of exactly the minimal elements $\rho_i$ of
each interval $j_i$ for $1\le i \le r$.
Finally, if $j_{r+1}\ne \emptyset $ and $\tau h = \infty $, then 
match $h \sim h\setminus \{ \rho_{\infty} \} $ for each face
$h\in e^{-1}(F_k)$ containing the element $\rho_{r+1 } $, and
match $h\sim h\cup \{ \rho_{r+1} \} $ otherwise. 

\begin{thm}\label{filter_match}
The above construction 
gives an acyclic matching on $\Delta (P)$.
\end{thm}
\proof
To check that our matching is acyclic, we 
verify (1) that it
is acyclic on each $F_k \setminus (\cup_{i<k} F_i )$ and 
(2) that there are no directed cycles involving multiple pieces of our 
decomposition, i.e. multiple fibres.  Fact (2) follows from the observation
that every edge between distinct fibres is oriented from the later fibre
to the earlier one.  Thus, we may henceforth restrict attention to a
single fibre $e^{-1}(F_k)$.

If the $I$-intervals for $F_k$ do not cover all ranks, then 
(1) follows from the fact that each upward step in a cycle must insert 
$\rho_0 $ while no downward step may delete $\rho_0 $.  It remains to 
prove (1) when $j_0 = \emptyset $.  The matching then breaks into two 
stages: (a) the matching on cells which do not include just
the minimal element for one or more $J$-intervals, and (b) matching on 
$2^u$ remaining cells which do include exactly the minimal element of 
each $J$-interval, letting $u$ be the number of ranks not
covered by any $J$-interval.  

In (a),
each face is matched with at most one other face by the bijection
between what we will call ``hollow'' cells and ``doubly-hit'' cells within
any fibre $e^{-1}(F_k)$.  A cell is {\bf hollow} if it lacks the minimal 
rank vertex in the first $J$-interval where it differs from the critical cell;
a cell is {\bf doubly-hit} if instead it includes the minimal rank
vertex and at least one other vertex in this $J$-interval.  
Our observation above that every face in
$F_k\setminus (\cup_{i<k} F_i )$ hits the lowest $J$-interval where it
differs from the potential
critical cell guarantees that this indeed gives a matching
on the faces in $F_k\setminus (\cup_{i<k} F_i )$.

First we show there are no directed 
cycles involving only cells matched in fashion (a).
Suppose there is such a cycle,
and consider any upward oriented edge $a\prec b$ in the cycle.
There must then be a downward edge from $b$ to some other $c$.  
A cycle may only
involve two ranks of the face poset since two consecutive upward edges
in a cycle are 
precluded by the up edges coming from a matching.  Thus, $a$ and $c$
must be hollow while $b$ must
be doubly-hit.  Furthermore, $c$ is obtained from $b$ by deleting some rank
above the rank where $a$ and $b$ differ or else $c$ would belong to an 
earlier orbit.  This choice of deletion precludes $c$ being hollow unless
we obtain $c$ by deleting from $b$ one element on the interval where $a$ and 
$b$ first differ, and $c$ is left with only the minimal element on that 
interval.  Continuing in this fashion, our cycle must produce cells that 
agree with the critical cell (if there is one) to higher and higher ranks, 
and thus the cycle may never be completed.  The point is that we may never 
move downward in the face poset by 
deleting the unique element covering an interval and also 
may never move upward by adding to an interval where we only have the minimal
element, so such intervals may never be altered.

On the other hand, there are no directed cycles 
involving only matching edges of
type (b), since the complete matching on a Boolean algebra $B_n$ which pairs
$S\setminus \{ i\} $ with $S\cup \{ i\} $ for some fixed $i\in [n]$ is
acyclic.  Finally, there are no directed cycles involving both types
of cells because there are no directed edges from cells of type (b)
to cells of type (a).  

\EOP 

Each fibre $e^{-1}(F_k)$ will have at most one critical cell, 
as described next.

\begin{thm}
Each facet contributes at most one critical cell.  The set of faces 
in $F_k \setminus (\cup_{i<k} F_i)$ includes 
a critical cell exactly when the attachment
of $F_k$ changes the homotopy type of the complex, i.e. when every vertex
in $F_k$ belongs to at least one of the intervals in $J(F_k)$.
\end{thm}
\proof
If some node is not in any minimal skipped interval of $F_k$, then the 
minimal rank such node $v$ is a cone point in 
$F_k \cap (\cup_{i<k} F_i)$ and 
is used to match faces including $v$
with those excluding $v$, so $F_k \setminus (\cup_{i<k} F_i )$ has
no critical cells.  If every node is in a minimal skipped interval of
$F_k$ but some node is not in any $J$-interval, again we obtain a 
complete acyclic matching on $e^{-1}(F_k)$.  In either case, the
complete acyclic matching on $e^{-1}(F_k)$ implies that  
the complex $F_1 \cup \cdots \cup F_k$
of faces belonging to one of the first $k$ facets collapses onto 
$F_1\cup \cdots \cup F_{k-1}$, so the homotopy type is unchanged.
If every node in $F_k$ belongs to at least one of its 
$J$-intervals, then our acyclic matching has a single 
critical cell comprised of the minimal
ranks of the $J$-intervals.  Thus, exactly one Morse number changes,
implying that the reduced Euler characteristic changes and hence the 
homotopy type also must change.  Alternatively, one may apply
the next proposition.
\EOP

\begin{cor}
When $F_j$ contributes a critical cell, then the dimension of the 
critical cell is one less than the number of $J$-intervals for $F_j$.
\end{cor}

The following proposition follows from a similar argument
to our proof of Theorem ~\ref{filter_match} (though this was not the
approach of Billera and Hetyei).  That is, one may 
give an acyclic matching on each $F_j \cap (\cup_{i<j} F_i)$ such that
(1) $F_j \cap (\cup_{i<j} F_i)$ has no critical cells when 
$F_j \setminus (\cup_{i<j} F_i )$ has no critical cells, and (2) 
otherwise, the critical cells in $F_j \cap (\cup_{i<j} F_i)$ are 
the proper faces of the unique critical cell in $F_j \setminus 
(\cup_{i<j} F_i)$.  

\begin{prop}[Billera-Hetyei]\label{overlap}
Given a lexicographic order $F_1,\dots ,F_k$, the 
intersection $F_j\cap (\cup_{i<j} F_i)$ is either homotopy equivalent 
to a sphere or is collapsible.
\end{prop}

From this, a simple Mayer-Vietorus argument yields:

\begin{cor}
Each facet attachment changes at most one Betti number, and changes
its value by at most one.
\end{cor}

We conclude this section with a useful fact for proving the uniqueness of 
a gradient path from a critical cell $\tau^{(p+1)}$ to a critical cell 
$\sigma^{(p)}$.  
First, it is helpful to note 
that any gradient path from $\tau^{(p+1)}$ to $\sigma^{(p)}$
must alternate between only
the ranks $p+1$ and $p$ (cf. [Fo]) since it cannot have two consecutive 
upward steps (since the upward edges comprise a matching) and also 
cannot end with an upward step (since $\sigma^{(p)} $ is not in the matching,
so there are no upward edges leading to $\sigma^{(p)} $).

\begin{prop}
If critical cells $\sigma^{(p)} $ and $\tau^{(p+1)}$ agree up through 
rank $r$ 
then every gradient path from $\tau $ to $\sigma $ leaves 
ranks $1,2,\dots ,r$ unchanged.
\end{prop}

\proof
If a rank $r'$ below $r+1$ is deleted in a gradient path
step $\tau_i \rightarrow \sigma_i$, then $e\sigma_i <_{lex} e\sigma $. 
This is because the first gradient path step
deleting an element at any rank $r' \le r$ will be immediately
followed by an upward step inserting an element at 
rank $r'$ or lower that belongs to a 
strictly earlier saturated chain.  To see this, recall that the element
at rank $r'$ was the lowest element of $e\tau $ in the interior of one
of its $J$-intervals, and that none of the earlier upward (matching) steps
will ever insert another element on this $J$-interval.
Thus, the chain obtained by deleting rank $r'$
belongs to a saturated chain which is strictly earlier in lexicographic
order than both $e\tau $ and $e\sigma $.
This makes it impossible for our Morse function to decrease at every 
remaining step on the gradient path to $\sigma $, a contradiction.
\EOP

In Section 5, we shall establish a notion of rank inversion 
for the poset of partitions
of a multiset in terms of labels on its covering relations.  

\section{M\"obius functions}
Let us observe, via the Morse inequalities and our description of 
critical cells in terms of skipped intervals in saturated chains, 
that the M\"obius 
function $\mu (\hat{0},\hat{1})$ may be computed by the following formula:
\begin{prop}
$$\mu_P(\hat{0},\hat{1}) = \sum_{\sigma \in Crit(\Delta )}
(-1)^{|J(e\sigma )| + 1},$$ 
where the sum is over critical cells $\sigma $ of a 
lexicographic discrete Morse
function, letting $J(C)$ denote the collection of 
disjoint skipped intervals in a saturated chain
$C$ obtained from the original skipped intervals $I(C)$.
\end{prop}

\proof
Consider the alternating sum 
$\sum_{i=-1}^{dim (\Delta)} (-1)^i m_i$
of criticals cells of various dimensions in a reduced discrete Morse
function, i.e. one which includes the empty set.
Now recall from the Morse inequalities:
$$\sum_{i=-1}^{dim (\Delta )} (-1)^i m_i = \sum_{i=-1}^{dim (\Delta )}
(-1)^i b_i = \tilde{\chi }(\Delta ).$$
Each saturated chain $m$ has at most one critical cell $\sigma $, and
this satisfies 
$\dim (\sigma ) = |J(m)| - 1$.  
\EOP

The facets $m$ contributing
critical cells are exactly those where every vertex in $m$ belongs to some
interval in $J(m)$.  Thus, one might compute $\mu_P(\hat{0},\hat{1})$ by 
characterizing facets which have a critical cell and determining 
dimension parities.  This may be an easy way (for instance) to show
the M\"obius function is 0 for certain classes of poset intervals; we give
one example of such a computation for the poset of partitions of a 
multiset later.

\begin{rk}
To obtain $\mu_P(u,v)$, simply
restrict the lexicographic order to the interval $(u,v)$ and consider
the resulting lexicographic discrete Morse function there.
\end{rk}  

One could also
derive a lexicographic order on any rank-selected subposet $P^S $, which in
turn would give a lexicographic discrete Morse function on $P^S $, giving an
expression for the rank-selected M\"obius function.
We have not looked into whether there is any useful 
relationship between
the critical cells on intervals or rank-selected subposets
and the critical cells for the entire poset.

\section{Relation to shellability}\label{shell_section}

This section makes explicit how to construct a discrete Morse function
from an arbitrary shelling; lexicographic discrete
Morse functions generalize the notion of a lexicographic 
shelling in exactly the same fashion.

\begin{defn}
A simplicial complex $\Delta $ is {\bf pure} of dimension $d$ if all its
maximal faces are $d$-dimensional.  A {\bf shelling} is a facet 
ordering $F_1,\dots ,F_k$ on $\Delta $ such that $F_j \cap (\cup_{i<j} F_i)$
is pure of codimension one for $2\le j \le k$.  
\end{defn}

A shelling gives a way of building up a pure complex $\Delta $ by 
sequentially inserting facets such
that each facet either attaches along its entire boundary or leaves the 
homotopy type unchanged, implying that $\Delta $ is collapsible or 
homotopy equivalent to a wedge of spheres of top dimension (and has a 
similar local property called Cohen-Macaulayness).  The following is very 
similar to [Ch, Theorem 4.2], but is included for 
completeness sake.

\begin{prop}
If $F_1,\dots ,F_k$ is a shelling order for the facets in a pure 
simplicial complex $\Delta $ of dimension $d$, 
then $\Delta $ has a discrete Morse function whose critical cells are
the facets $F_j$ of $\Delta $ such that $F_j \cap (\cup_{i<j} F_i )$ is
the entire boundary of $F_j$.
\end{prop}
\proof
Since $F_1,\dots ,F_k$ is a shelling, $F_j \cap (\cup_{i<j} F_i)$ is pure
of codimension one, so its maximal faces each omit a single vertex of $F_j$.
Hence, the intersection is either the entire boundary of $F_j$ or it has
at least one cone point, namely a vertex which is not omitted by any 
of the maximal faces in $F_j \cap (\cup_{i<j} F_i)$.  In either case, we 
give a matching on the faces in $F_j \setminus (\cup_{i<j} F_i)$ so that the 
union of these matchings is an acyclic matching on $\Delta $ with the desired
set of critical cells unmatched.  
When $F_j \cap (\cup_{i<j} F_i)$ has one or more 
cone points, choose any such cone point $v$.  Match
each face in $F_j\setminus (\cup_{i<j} F_i)$ which includes $v$ with the face
in $F_j \setminus (\cup_{i<j} F_i)$ obtained by deleting $v$.  
Note that every
face in $F_j\setminus (\cup_{i<j} F_i)$ is used in this matching 
exactly once.
When $F_j\cap (\cup_{i<j} F_i)$ is the entire boundary of $F_j$, then 
the only
face in $F_j\setminus (\cup_{i<j} F_i )$ is the facet $F_j$, and it is 
left unmatched.  

It is easy to check that the union of these matchings is
an acyclic matching on the face poset, since matching edges (i.e. the
edges which
get directed upward) are covering relations $u\prec v$ where $u,v$ are both
contributed by the same facet.  Hence, there 
are no covering relations $u\prec v$ with $v$ belonging to an 
earlier facet than $u$ does.  Hence, any directed cycle would live
in some $F_j \setminus (\cup_{i<j} F_i)$.  
This is impossible because each of the upward oriented edges in 
$F_j \setminus
(\cup_{i<j} F_i )$ inserts the same cone point $v$ but no downward 
edges may delete
$v$.
\EOP

Whenever a lexicographic order on facets gives rise to a shelling, we call
this a lexicographic shelling.  This is equivalent to Kozlov's notions of
EC-shellability (for edge-labelings) and CC-shellability (for 
chain-labelings), as introduced in [Ko1].  

\begin{prop}
A lexicographic order on facets is a shelling order if it gives 
rise to a lexicographic discrete Morse 
function in which all minimal skipped intervals have height one.
\end{prop}
\proof
Specialize the above proof to this situation to get a lexicographic discrete
Morse function.  The reverse direction is shown in [Ko1].
\EOP

Kozlov's notion amounts to requiring for each interval that
all but the lexicographically earliest chain have a minimal skipped interval
of height one.  Such an interval of height one is called a 
topological descent in [HK], [He1] because 
it behaves topologically like a descent would behave in an
EL-labeling or CL-labeling.

\section{On the homotopy type of $\Delta 
(\Pi_n/S_{\lambda })$}\label{homotopy}

This section uses a particular 
lexicographic discrete Morse function to show for 
$\lambda $ of hook shape that the poset 
$\Delta(\Pi_n/S_{\lambda })$ of partitions of a multiset ordered by
refinement is either collapsible or homotopy equivalent to
a wedge of spheres of dimension $n-3$.  It seems likely
that with further work this hook shape assumption may be relaxed or perhaps
even removed. 

\begin{figure}[h]
\begin{picture}(250,155)(0,10)
\includegraphics[width=0.7\textwidth, angle = 0]{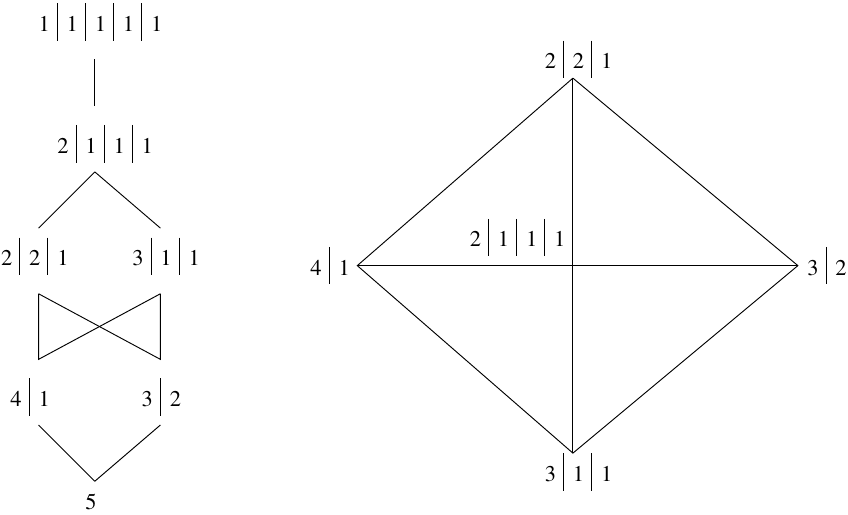}
\end{picture}
\caption{$\Pi_5/S_5 $ and its (2-dimensional) order complex $\Delta ( \Pi_5/S_5)$}
\label{pi5mod_s5}
\end{figure}

Let $\Pi_n$ denote the poset of set partitions of $\{ 1,2,\dots ,n\} $
ordered by refinement.  Let $\lambda = (\lambda_1 ,\dots ,\lambda_k )$ be 
a number partition of the integer $n$ into unordered 
parts $\lambda_1 ,\dots ,\lambda_k $.  
By convention, we list the parts in decreasing order, so they satisfy
$\lambda_1 \ge \cdots \ge \lambda_k > 0$ and $\lambda_1 + \cdots + 
\lambda_k = n$.  Denote by $S_{\lambda }$ the Young subgroup 
$S_{\lambda_1 }\times \cdots \times S_{\lambda_k}$ of the symmetric group
$S_n$.  
Notice that $S_{\lambda }$ acts on $\{ 1,2,\dots ,n\} $ by letting
each $S_{\lambda_i }$ permute the values $\{ \lambda_1 + \cdots +
\lambda_{i-1} + 1, \dots ,\lambda_1 + \cdots + \lambda_{i-1} + \lambda_i \}$
among themselves.  This naturally gives rise to an order-preserving
$S_{\lambda }$-action
on the set partitions of $\{ 1,2,\dots ,n\} $, i.e. an action such 
that $u<v$ implies $\sigma u < \sigma v$ for each $\sigma\in S_{\lambda }$
and $u,v$ comparable elements of $\Pi_n$.  Thus, it makes sense to consider
the quotient poset $\Pi_n /S_{\lambda }$, i.e. the poset of partitions of
a set with repetition  
$\{ 1^{\lambda_1},\dots ,k^{\lambda_k } \} $ (which we call a 
{\bf multiset}).  
\begin{figure}[h]
\begin{picture}(250,155)(0,10)
\includegraphics[width=0.75\textwidth, angle = 0]{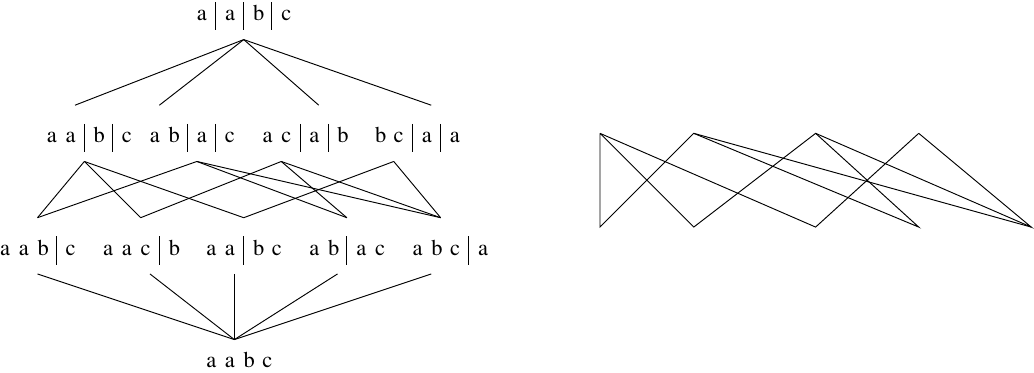}
\end{picture}
\caption{$\Pi_4/S_2\times S_1\times S_1$ and its (1-dimensional) order complex}
\label{pi4mod_s211}
\end{figure}
Figures ~\ref{pi4mod_s211} and ~\ref{pi5mod_s5} 
show the posets $\Pi_4 /S_2\times S_1\times S_1$ 
and $P_5 \cong \Pi_5/S_5$, respectively, along with their order complexes.

We will study $\Pi_n /S_{\lambda }$ in the special
case where $\lambda $ is a number partition of {\bf hook shape}, namely 
one in which $\lambda_2 = \cdots = \lambda_k = 1$.  Notice that this
extrapolates between the partition lattice $\Pi_n$ (in which case $\lambda_1$
also equals 1) and the poset $P_n$ of number partitions of an integer $n$
ordered by refinement (in which case $k=1$).  

Ziegler showed that $\Delta(\Pi_n/S_n)$ is not Cohen-Macaulay for $n\ge 19$
by exhibiting a disconnected interval of rank 3 in the poset $P_{19}$ 
[Zi, p. 218].  
Recall that $\Delta $ is {\bf Cohen-Macaulay} if the link of each face
(including the empty face) only
 has top dimensional (reduced) homology; note that a disconnected graph
$G$ has dimension 1, but $\tilde{H_0}(G) \ne 0$, so it is not Cohen-Macaulay.
At the other extreme, Gessel gave an EL-labeling 
for $\Pi_n $, implying $\Delta (\Pi_n )$ is homotopy
equivalent to a wedge of spheres of top dimension, as discussed in
[Bj1, p. 165].  We raise 
the question of whether all posets $\Pi_n /S_{\lambda }$ are nonetheless
either collapsible or homotopy equivalent to a wedge of spheres of top
dimension, and give an affirmative answer for $\lambda $ of hook shape.

\begin{examp}[Ziegler]\label{zieg}
The open interval $(6|5|3|2|2|1,
8|7|4)$ prevents $P_{19}$ from being 
Cohen-Macaulay.  Simply note that 
$8|7|4$ may be refined to
$6|5|3|2|1|1$ in two different ways:
\begin{alignat*}{2}
8 &= 6+2 &= 5+3\\
7 &= 5+2 &= 6+1\\
4 &= 3+1 &= 2+2
\end{alignat*}
Any saturated chain
from $\hat{0} $ to $6|5|3|2|2|1$ together with any
saturated chain from $8|7|4 $ to $\hat{1} $ gives a face 
$F\in \Delta(P_{19})$
such that $lk(F)$ is a disconnected graph, and hence has lower homology.
\end{examp}

The remainder of this paper will be concerned with proving the following
theorem and deducing a consequence about M\"obius functions.

\begin{thm}\label{hooksphere}
If $\lambda $ is a partition of hook-shape then the poset $\Pi_n/S_{\lambda }$
of partitions of the multiset $\{ 1^{\lambda_1},\dots ,k^{\lambda_k }\} $
ordered by refinement is either collapsible or
has the homotopy type of a wedge of spheres 
concentrated in top dimension.
\end{thm}

Here is an outline of the proof:

\begin{enumerate}
\item
Give a {\bf lexicographic order} on saturated chains based on a 
chain-labeling.  We shall then consider the resulting lexicographic
discrete Morse function.
\item
{\bf Classify} lower critical cells (i.e. critical cells
not of top dimension) by classifying the three ways in which 
nontrivial minimal skipped intervals may arise in posets
$\Pi_n/S_{\lambda }$.  This will enable us to 
pair lower critical cells based on the structure of their highest rank
nontrivial minimal skipped interval.  One of these classes of nontrivial
minimal skipped interval cannot occur in critical cells 
for $\lambda $ of hook shape.
\item
For the two types of lower critical cells $\sigma $ that occur in hook
shapes, we prove (by explicit construction) the {\bf existence of a partner 
critical cell} $\tau $ such that either there is a gradient
path from $\sigma $ to $\tau $ (for $\dim (\sigma ) = \dim (\tau ) +1$) 
or from $\tau $ to $\sigma $ (for $\dim (\tau ) = \dim (\sigma ) +1$).  
For the remaining class of lower critical cells, we still
construct partner cells, and it is open 
whether or not these partner cells are always critical.
\item
Prove that for each pair $\sigma ,\tau $ of partner (critical) cells 
there is a {\bf unique gradient path} from $\tau $ to
$\sigma $.
This implies that $\sigma $ and $\tau $ may indeed cancel by reversing this
gradient path to obtain a larger acyclic matching.
This uniqueness of gradient paths follows from our use of a statistic 
called the {\bf inversion set} of a cell which 
decreases at every step in a gradient path.  
\item
Choose how to {\bf pair critical cells} (in a way related to step 4) so that
no two critical cell have the same partner.
Furthermore, we show that partner critical cells occur consecutively 
among the set of all critical cells in lexicographic order, so cancellation
proceeds as desired.
\end{enumerate}

The third step
in the proof uses the assumption that $\lambda $ is of hook-shape, so 
at the end of that section we indicate what exactly would be
needed to generalize the proof to all $\lambda $.  It seems likely that the 
hook-shape assumption can at least be weakened.  
Now let us turn to the details of the proof.

\subsection{A lexicographic order on saturated chains}

We will use a chain-labeling to order the saturated chains, so each covering
relation $v_{i-1}\prec v_i$ is given a label which may depend on the 
rooted chain $\hat{0}\prec v_1\prec\cdots\prec v_{i-1}$ which it extends to
$\hat{0}\prec\cdots
\prec v_{i-1}\prec v_i$; saturated chains are ordered by the lexicographic
order on their label sequences $\lambda_1\lambda_2\cdots\lambda_n $ where 
$\lambda_i$ is the label for $v_{i-1}\prec v_i$.  In particular, the 
label we assign to $u\prec v$ in $\Pi_n/S_{\lambda }$ will depend on an
ordering of the blocks of $u$ which is determined inductively by the choice
of root for $u$. That is, we obtain a block order for $v_i$ from a block order
for $v_{i-1}$ whenever $v_{i-1}\prec v_i$ by preserving the order in $v_{i-1}$
in all but the one block which gets refined to obtain $v_i$ and replacing this
block by the two resulting blocks as follows: if there is more than one block
in $v_{i-1}$ of the type to be refined then we refine the leftmost such block,
and when replacing it by the two new blocks we place the smaller one to the 
left of the larger one.  We determine which of these two blocks is 
smaller by using the {\bf length-lex} block order (defined in [HK]), namely by
letting $B_1 < B_2$ for $|B_1| < |B_2|$ or for $|B_1|=|B_2|$ and $w_{B_1} 
< w_{B_2}$ in the lexicographic order on words
where $w_B$ is the word obtained from the 
multiset $B$ by arranging the letters in increasing order.  

When a covering relation $u\prec v$ splits a block $B$ into children 
$B_1,B_2$ such that $B_1 < B_2$ in length-lex order, let us denote $B_1$ 
by $L_{u\prec v}$ and $B_2$ by $R_{u\prec v}$, since $B_1$ is the left
child and $B_2$ is the right child.  
The rooted covering relation $v_{i-1} \prec v_i$ is then labeled by the 
ordered pair $(n_i,L_{v_{i-1}\prec v_i})$ 
consisting of the position $n_i$ of the 
newly inserted bar together with the block word 
immediately to the left of the new bar in $v_i$; 
the bar position $n_i$ has higher precedence in the label.
As an example of our labeling, the saturated chain
$$aabbb < ab|abb < ab|b|ab < a|b|b|ab < a|b|b|a|b $$ in $\Pi_5/S_{(3,2)}$
is labeled $((2,ab),(3,b),(1,a),(4,a))$.  We will often use a more 
compact notation for saturated chains, namely we will write the most 
refined element of the chain with bars separating the elements labeled 
by the rank at which the bars are inserted in the course of successively
refining the partition.
For instance, the above chain could be more succinctly expressed as
$a|_3 b|_1 b|_2 a|_4 b$.

\subsection{Classification of lower critical cells in terms of 
skipped intervals}

When a critical cell is not a facet, it must have at least one nontrivial 
minimal skipped interval (i.e. one that skips multiple ranks).  In this 
section, we classify the critical cells which are of smaller than top
possible dimension according to the structure of their nontrivial minimal
skipped interval appearing at the highest rank.  We will 
refer to critical cells which 
are not top-dimensional as {\bf lower critical cells}.

There are three fundamental ways in which nontrivial minimal skipped 
intervals may arise (as well as hybrid mixtures of these).  An instance
of type three is discussed in example ~\ref{critpair} in later sections.
In each case, we have a saturated chain $u=u_1 \prec u_2 \prec\cdots\prec
u_k = v$ on an interval $(u,v)$ which is not lexicographically smallest on
$(u,v)$ but where $u_i\prec u_{i+1} \prec u_{i+1} $ is lexicographically 
smallest on the interval $(u_i,u_{i+2})$ for each $1\le i< k-2$.  It will be 
useful for us to begin by describing the two ways that extensions
$u_i\prec u_{i+1}
\prec u_{i+1}$ of the interval $(u_i,u_{i+2})$ might 
not be lexicographically smallest on the interval 
$(u_i,u_{i+2})$.

\begin{enumerate}
\item
when bars are inserted sequentially from right to left in consecutive
steps $u_i\prec u_{i+1}$ and $u_{i+1}\prec u_{i+2}$.
\item
when a single block of $u_i$ is split into three pieces in $u_{i+2}$ 
such that bars are inserted from left to right in the chain $u_i\prec u_{i+1}
\prec u_{i+2}$ with the resulting blocks decreasing in size from left to
right, i.e. with $L_{u_i\prec u_{i+1}} > L_{u_{i+1}\prec u_{i+2}}$.
\end{enumerate}

Following [HK], we say there 
is a {\bf topological descent} at $u_{i+1}$ in these 
situations where $u_i\prec 
u_{i+1}\prec u_{i+2}$ is not lexicographically smallest on $(u_i,u_{i+2})$.
As examples, 
the chain $a|_2a|_1aa$ has a topological descent of type 1 
at rank 1 while the chain $b|_1a|_2c$ has a topological descent of type 2
at rank 1.  
Since every rank in $e\sigma $ must belong to some interval in 
$I(e\sigma )$ for $\sigma $ to be critical, there must be a topological
descent at every rank above the highest rank nontrivial interval in $I$.

Let us now list the three basic forms which nontrivial minimal
skipped intervals $(u,v)$ may take.
To get such an interval $(u,v)$, there 
must be saturated chains $C_1, C_2$ both containing $u,v$
which agree outside of $(u,v)$ where $C_1$ is lexicographically 
smaller than $C_2$ on $(u,v)$ but where neither chain has any topological
descents on $(u,v)$; let $C_2$ include $u=u_1\prec\cdots
\prec u_k=v$.  Then $C_1$ and $C_2$ must be related
to each other in one of the following ways:
\begin{enumerate}
\item
The chains $C_1$ and $C_2$ both refine a pair of 
identical blocks $B_1,B_2$ of $u$ 
into pieces $\alpha_1,\dots ,\alpha_r, \beta_1,\beta_2$ in $v$ such
that
\begin{enumerate}
\item
$\beta_1,\beta_2$ are the children of $B_i$ in $C_i$ with $\beta_1 $
smaller than $\beta_2 $ in our block order
\item
$\alpha_1,\dots ,\alpha_r$ are the children of $B_1$ in $C_2$ and of 
$B_2$ in $C_1$ with $\alpha_1,\dots ,\alpha_r$ arranged in increasing
order in our block order
\item
$r > 2$
\item
$\alpha_1 $ is larger than $\beta_1$ in our block order
\end{enumerate}
The blocks
$B_1$ and $B_2$ are distinguished from each other by the fact that 
$B_1$ is refined before $B_2$ because it is farther to the left in the ordered 
partition $u$ (where the ordering of blocks in $u$ results from the choice
of root that is common to $C_1$ and $C_2$).  
Both $C_1$ and $C_2$ refine the blocks of $u$ from 
left to right with the children increasing in length-lex order from left to 
right.  Figure ~\ref{skipint1} 
gives an example with $u$ of rank 1, $v$ of rank 4, 
and the $C_1$ refinement represented by the top arrows 
and the $C_2$ refinement by the bottom arrows, both denoting ranks when 
bars are inserted in these positions.
\begin{figure}[h]
\begin{picture}(250,70)(0,0)
\includegraphics[width=0.7\textwidth, angle = 0]{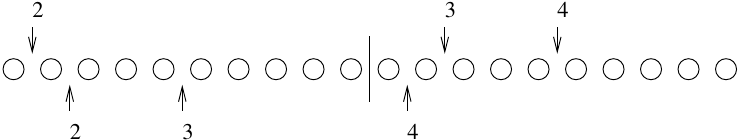}
\end{picture}
\caption{An example of type 1}
\label{skipint1}
\end{figure}
\item
The chains $C_1$ and $C_2$ both
split the same set of blocks of $u$ into the same set of smaller blocks in
$v$ in two different ways by rearranging which 
children blocks belong to which types of parent.  Ziegler constructed an
interval of this type to show that the poset of partitions of an 
integer is not Cohen-Macaulay [Zi]; he noted that $4=1+3=2+2;
7=2+5=1+6; 8=2+6=3+5$, allowing the blocks $4,7,8$ to be split into the 
six pieces $1,3,2,2,5,6$ in two different ways, yielding two disjoint 
boolean algebras (and consequently 
two different saturated chains on the interval $(4|7|8,2|2|1|6|3|5)$
both of which are free of topological descents), as depicted in
Figure ~\ref{skipint2}. 
Again we list $C_1$ above $C_2$ in the diagram.  
\begin{figure}[h]
\begin{picture}(250,70)(0,0)
\includegraphics[width=0.7\textwidth, angle = 0]{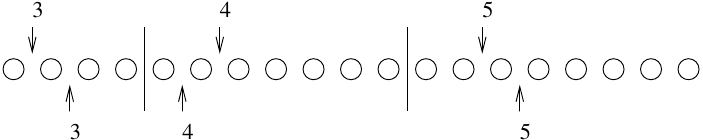}
\end{picture}
\caption{An example of type 2}
\label{skipint2}
\end{figure}
\item
Some block $B$ gets refined in $C_1$ but not in $C_2$ (or the 
reverse), because in the refinement in $C_1$, the block 
$B$ is a child of some larger block $B'$ that also appears
in $u$ and so the copy of $B$ in $u$ is refined and a new copy is 
recovered
as a child of $B'$ thereby appearing both in $u$ and $v$.  
Figure ~\ref{skipint3} gives an example for $u$ of rank 2, $v$ of rank 5
with $B=1^2, B'=1^{10}$.
\begin{figure}[h]
\begin{picture}(250,70)(0,0)
\includegraphics[width=0.7\textwidth, angle = 0]{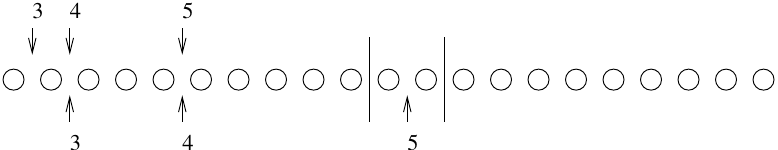}
\end{picture}
\caption{An example of type 3}
\label{skipint3}
\end{figure}
\end{enumerate}

It is not hard to see that these are the only basic possibilities given that
$C_2$ cannot have any topological descents, so 
bars must be inserted left to right within $u$ with the children of each 
$u$-block increasing in length-lex order from left to right to obtain $v$.
Case one is the only way the same types of blocks may
be split into the same types of pieces, i.e. by rearranging which actual
blocks in $u$ of a particular content are split which ways.  Case two
handles the possibility that the same set of blocks of $u$ gets split into
the same set of smaller blocks in $v$, but with different types of blocks
having different types of parents in $C_1$ and $C_2$.
Then case three is the only way different actual
blocks of $u$ could be refined with $C_1$ and $C_2$ both leading to the 
same set of blocks in $v$.  Case two cannot happen for $\lambda $ of hook 
shape because it requires repetition of more than one type of letter in order
for the same child blocks to come from more than one possible parent.

Every lower critical cell in a lexicographic discrete Morse function 
must have at least one nontrivial minimal skipped
interval in $J$, and we find it convenient to pair lower critical cells
based on their highest rank such interval.
In each of the three cases, we will construct a partner cell by showing that
for any particular highest rank nontrivial minimal skipped interval
$(u,v)$ exactly one of the following two possibilities holds:
\begin{itemize}
\item
there is a bar inserted to refine $v_i$ to get $v_{i+1}$ in 
$C_2$ within the interval $(u,v)$ whose insertion
may be shifted to a rank coming after the interval $(u,v)$, 
yielding a lexicographically later saturated chain $C_3$ 
which has a skipped interval of height one smaller than
$(u,v)$ and which has only topological descents after this skipped
interval.
\item
the interval $(u,v)$ can be enlarged
by shifting a single bar insertion into the 
interval from a later rank, yielding a lexicographically earlier
saturated chain $C_2'$ which has only topological descents
after the interval that replaces $(u,v)$, 
so this has a (potentially critical) 
cell of dimension one smaller than the critical cell in $C_2$.
\end{itemize}

Either possibility 
guarantees the existence of a partner cell for each lower critical 
cell because the lower ranks of $C_2$ are unchanged, so that minimal 
skipped intervals in $I(C_3)$ (resp. $I(C_2')$) cover $C_3$ (resp. $C_2'$).
We will also show in a later section that there is a unique gradient
path between these partner cells, as needed for cancellation.  For the
nontrivial minimal skipped intervals of types one and three, we will show
that these partner cells are themselves critical, but it is not known
whether this is always true
for nontrivial minimal skipped intervals of type two.

\subsection{Existence of partner critical cells}\label{partner}

For each of the three types of nontrivial minimal 
skipped interval which may occur in $\Pi_n/S_{\lambda }$
as the highest rank such interval in some $e\sigma $, we describe
how a partner cell $\tau $ is constructed and for types 1 and 3 we observe
that this partner must always be critical.  
The constructions of partner cells 
generalize in an obvious way to mixtures of these ``fundamental'' types of 
nontrivial minimal skipped intervals, and mixtures of types 1 and 3 have
critical partners for the same reasons that cells with intervals purely of
type 1 or 3 have critical partners.  

\begin{figure}[h]
\begin{picture}(250,100)(-30,0)
\includegraphics[width=0.5\textwidth, angle = 0]{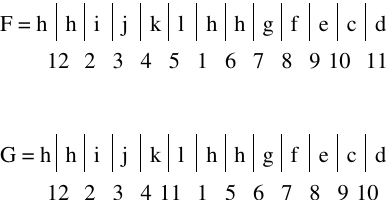}
\end{picture}
\caption{Two saturated chains whose critical cells are partners}
\label{pairing}
\end{figure}

Let us begin with an example (where the nontrivial minimal skipped interval
is of type 3), depicted in Figure ~\ref{pairing}.

\begin{examp}\label{critpair}
\rm
The saturated chain $F$ has nontrivial minimal skipped intervals with rank
sets $\{ 1\}, \{ 2\} , \{ 3, 4, 5, 6\} ,\{7\} , \{ 8\} , \{ 9\} , \{ 10\} 
, \{ 11\} $ yielding a critical cell $\sigma^{(6)} $ consisting of ranks 
$\{ 1, 2, 3, 7, 8, 9, 10\} $.  
The unique nontrivial minimal 
skipped interval in $F$ consists of rank set $\{ 3, 4, 5, 6\} $ and results
from the block $\{ h, h\} $ also being contained in the block $\{ h, h, g, f,
e, c, d\} $ so that a lexicographically earlier saturated chain may
sequentially insert bars into positions $6, 2, 1, 3, 4, 5, 8, 9, 10, 11, 7$
splitting the block $\{ h, h\} $ into singletons immediately after its 
creation instead of later splitting off
the two copies of $h$ which occur farther to the right.
The lexicographically later saturated chain
$G$ has nontrivial minimal skipped intervals with rank sets
$\{ 1\}, \{ 2\} , \{ 3, 4, 5 \} ,\{ 6\} , \{7\} , \{ 8\} , \{ 9\} , \{ 10\} 
, \{ 11\} $  leading to a critical cell $\tau^{(7)} $ comprised of ranks 
$\{ 1, 2, 3, 6, 7, 8, 9, 10, 11\} $.  The gradient path from $\tau $ to 
$\sigma $ begins by deleting the descent at rank 10 in $G$ and replacing it
by an ascent; the gradient path continues by sequentially replacing ranks
$9, 8, 7$ and $6$ by ascents with the net effect of shifting the insertion 
of the bar between $k$ and $l$ from rank 11 to rank 5, increasing the height
of the nontrivial minimal skipped interval by one and thereby giving a new
critical cell of dimension one smaller.  
\end{examp}

We will always pair critical cells by
shifting bars in or out of the last nontrivial minimal skipped interval.  
A later section will use permutation statistics to 
prove uniqueness of gradient paths, at which time we will revisit
this example.  

The first type of highest rank
nontrivial minimal skipped interval for $e\sigma $ 
results from two identical blocks $B_1,B_2$ 
being refined in different ways, in which case $e\sigma $
refines $B_1$ in a lexicographically larger way than $\gamma $ 
does (letting $\gamma $ denotes 
the lexicographically earliest saturated chain 
differing from $e\sigma $ only on this
interval).  Because the interval is nontrivial,
either $e\sigma $ inserts multiple bars left to right in $B_1$  
creating children from left to right which are nondecreasing before
refining $B_2$, or
there is at least one block to the right of $B_1$ and to the left of $B_2$ 
which also gets refined in the interval.
In either case, let us refer to the series of bar insertions that are 
strictly within the skipped interval as $b_1,\dots ,b_r$.
All higher ranks in $e\sigma $ must consist
entirely of topological descents or skipped intervals which have height one
after truncation.  We defer discussion of truncation to the end of
Theorem ~\ref{uniquepath}, and at this point note that it will not change
our partnership construction because the descent to be shifted in or out
of the interval
occurs outside the interval at a higher rank than all such truncated 
intervals.  Assuming there are only topological descents above the 
interval, the remainder of $e\sigma $
must entirely refine each of the remaining blocks with only topological
descents and proceed from right to left among the blocks to be refined.
We will have two quite different cases, depending on how the left child
$L_r$ resulting from the refinement step $b_r$ compares to the left block 
$L_{R_r}$ obtained from the right child $R_r$ of the refinement $b_r$
when the block $R_r$ is first refined at some later point in the chain.

If $L_r > L_{R_r}$, then we pair $\sigma $
with a face $\tau $ such that $e\tau $ has the $b_r$ bar 
insertion shifted to after the skipped interval, increasing the number of 
minimal skipped intervals by 1, so that 
$\dim (\tau ) = \dim (\sigma ) + 1$ and 
$e\sigma <_{lex } e\tau $.  The bar insertion $b_r$ 
is shifted to the unique position later in the 
saturated chain which ensures there are still only topological 
descents after the skipped interval.  Such a position exists precisely
because $L_r > L_{R_r}$.  One may check that $\tau $ is critical by 
observing
that the skipped interval derived from the former highest rank
nontrivial minimal skipped interval of $e\sigma $ is also minimal (
of type 1); moreover, the minimal skipped intervals of $\sigma $ after
truncation still cover all ranks, for the same reason this holds for
$\sigma $.

When $L_r \le L_{R_r}$, then we instead pair $\sigma $ with a 
critical cell $\rho $ such that $e\rho <_{lex } e\sigma $ and $\dim (\tau ) 
= \dim (\sigma )-1$.  This 
is accomplished by shifting the bar insertion which creates 
the block $L_{R_r}$ into the
minimal skipped interval increasing its height by one and decreasing the 
number of minimal skipped intervals by one; again the property of 
truncated intervals covering all ranks is preserved.  
Because exactly one of the 
conditions $L_r > L_{R_r}$ or $L_r\le L_{R_r}$ holds, 
we are ensured that there is
not a series of more than two critical cells arising from 
shifting more and more bar insertions in or 
out of the interval, so there is no issue of multiple cells needing the 
same partner in this manner.

Now let us turn to the possibility that the last nontrivial minimal 
skipped interval is of type two, in which case
the interval comes from refining a collection
of blocks in two different ways,
such as in Ziegler's example.  
When $\lambda $ is of hook shape, none of the
critical cells may be of type two, since the blocks may 
only use one type of letter with repetition, but every letter appearing
in the blocks to be refined must occur in at least two places;
it is impossible
to avoid topological ascents after the last nontrivial minimal skipped 
interval when using only the one type of letter which occurs with 
repetition.  Let us nonetheless briefly discuss how to match critical
cells $\tau $ whose last nontrivial minmimal skipped interval is of type
two, in case this could be helpful for generalizing our result
to all $\lambda $.  What remains open is whether the partner cell for
$\tau $ is itself critical.
We pair the critical cell $\tau $ with a lexicographically earlier cell
$\sigma $ by
replacing the (topological) descent immediately following the last nontrivial 
minimal skipped interval by an ascent in a way so that this rank is
incorporated into the last nontrivial minimal skipped interval.  
This is possible because
the last block refined in the last nontrivial minimal skipped interval
must have size at least three, because it is refineable in two 
different ways.  Thus, this block needs further refinement after the 
nontrivial minimal skipped interval, so we insert such a refinement 
step into the interior of the skipped interval.  

In the third case,
there is a block of $u$ which is refined but then rearises as a child 
block later in one of the facets, while this same block is not 
refined at all in the other facet.  This case proceeds identically to case one,
so we omit the argument.  
The point is again to compare $L_r$ to $L_{R_r}$ to 
decide whether the partner critical
cell should be lexicographically earlier or later.

\subsection{Collections of inversions: a decreasing statistic}

For each saturated chain, let us 
consider its restriction to those ranks  
above its last nontrivial minimal
skipped interval.  To each such restricted chain $v_r\prec\cdots\prec v_n$
where each $v_k$ has rank $k$ we associate 
a collection of rank inversions $(i,j)$ as follows:
we include $(i,j)$ among the set of inversion pairs whenever
either of the following conditions are met:
\begin{itemize}
\item
the bar inserted at rank $j$ (i.e. by the covering relation $v_{j-1} \prec
v_j $) is farther to the right than the bar inserted at rank $i$
\item
the covering relations between ranks $i,i+1,\dots ,j-1,j$ insert bars left
to right into a single block of $v_{i-1}$ in such a way that $L_i > L_j$.
\end{itemize}

Note that for the collection of 
inversions associated to a saturated chain to be exactly the inversions
of some permutation, the following two facts are necessary:
\begin{enumerate}
\item
$(i,j),(j,k)\in S$ implies $(i,k)\in S$
\item
$(i,k)\in S$ implies that $(i,j)\in S $ or $(j,k)\in S$ for each $i<j<k$
\end{enumerate}
The second condition implies that
whenever $(i,k)\in S$ then there is also some $(j,j+1)\in S$ for 
$i\le j < k$, 
enabling any permutation to be expressed as a product of adjacent 
transpositions.  These two axioms do not hold for saturated chains in 
their entirety, but
do hold when we restrict to the ranks above the last nontrivial minimal
skipped interval, i.e. the ranks which are allowed to 
vary in gradient paths
between pairs of critical cells to be cancelled.

\begin{rk}
In an edge-labeling in which consecutive labels increase exactly where there
are topological ascents, rank-inversions may be defined by
comparing the labels at the respective ranks.  
This is not the situation in our chain-labeling for
$\Pi_n/S_{\lambda }$.
\end{rk}

\subsection{Uniqueness of gradient paths}\label{unique}

In this section, we check for critical cells $\tau^{(p+1)},\sigma^{(p)}$
that were paired in Section ~\ref{partner} that the gradient
path from $\tau $ to $\sigma $ is unique.  
This makes use of the statistic of the 
previous subsection.  We return to the example of Section ~\ref{partner}
immediately after proving the following theorem.

\begin{thm}\label{uniquepath}
If $\sigma ,\tau $ are paired as in Section ~\ref{partner}, then
there is a unique gradient path from $\tau $ to $\sigma $.
\end{thm}

\proof
Let $\tau^{(p+1)}$ and $\sigma^{(p)}$ be partner critical cells with 
$e\tau >_{lex }e\sigma $.
Consider the gradient path from $\tau $ to 
$\sigma $ as given in section 4.3.  
If $e\tau $ and $e\sigma $ agree on ranks 
$1,2,\dots ,i$ 
then gradient paths from $\tau $ to $\sigma $ may never delete
elements of rank less than or equal to $i$.  This is because any
such deletion 
$\tau_i\rightarrow \sigma_i $ would leave a minimal skipped interval of 
$e\tau $ uncovered so that $e\sigma_i <_{lex } e\sigma $, implying 
$f(\sigma_i ) < f(\sigma )$.  This would make it impossible to complete
the gradient path to $\sigma $ with $f$ decreasing at each step.
Thus, we need only consider gradient path
steps deleting ranks above the last
nontrivial minimal skipped interval, since $\sigma $ and $\tau $ agree up 
through the first rank in this interval.

To this end, we shall use inversions to prove uniqueness.  Let us temporarily
assume that there are only topological descents after the last nontrivial 
skipped interval.  Later we will consider the possibility of nontrivial 
minimal skipped intervals which are reduced to height one by truncation.
Given a gradient
path from $\tau $ to $\sigma $, each step 
in a gradient path eliminates a single inversion since the downward 
step $\tau_i\rightarrow \sigma_i$ must uncover a topological descent 
yielding
$e\sigma_i <_{lex} e\tau_i$ so that
our matching forces $\tau_{i+1}$ to have this rank which we call $r$
reinserted with a lexicographically earlier topological ascent eliminating
a single rank inversion $(r,r+1)$.
Denote by $d_{i_1}\circ\cdots\circ d_{i_k} $
a gradient path which deletes rank $i_j$ at the step 
$\tau_j\rightarrow \sigma_j$ for each $j$.  The inversion set of $\tau $ is
then the inversion set of $\sigma $ together with the set of inversions 
found in the permutation $s_{i_1}\circ\cdots \circ s_{i_k}$ where $s_i$ 
denotes the adjacent transposition $(i, i+1)$ swapping $i$ and $i+1$.
Every other gradient path from $\tau $ to $\sigma $ must similarly 
give rise to another reduced expression for this same permutation since
every gradient path from $\tau $ to $\sigma $ must eliminate the same
set of inversions.

Furthermore, we need for the 
indices in the word $d_{i_1}\circ \cdots \circ d_{i_k}$ to be monotonically
decreasing so that a downward step in 
a gradient path always leads to a vertex that is matched with one above 
rather than below it; the 
cyclic shifts which arise from shifting a single bar insertion into the last 
nontrivial minimal skipped interval from a later rank meet this requirement.  
Note that the downward step $d_i$ eliminates an inversion at rank $i$ and 
the upward step via the matching replaces this with a lexicographically 
earlier ascent at rank $i$, so then applying  
$d_j$ immediately after $d_i$ for $j>i$ would lead downward to a face 
$\sigma_k$ which includes rank $i$ as an ascent, hence as a cone point 
in the intersection of $e\sigma_k $
with earlier facets, so that $\sigma_k$ is the top of a matching edge in
an earlier fibre, making it impossible for us to continue our gradient
path with an upward step from $\sigma_k $.  We have paired 
critical cells $\sigma^{(p)},\tau^{(p+1)}$ so 
that the permutation associated to the gradient path from $\tau $ to 
$\sigma $ is a cyclic shift 
$s_js_{j-1}\cdots s_{i+1}s_i$ with decreasing indices which has
a unique reduced expression, as desired.  

Now let us deal with the possibility of nontrivial skipped intervals in 
$I(e\sigma )$ or $I(e\tau )$ which become 
trivial in $J$ and which come after the last nontrivial
interval of $I$ that remains nontrivial in $J$.  (Let us call this last 
interval $\gamma $.)
Notice first that these must all occur consecutively immediately after 
$\gamma $, and after this we have only 
topological descents.  The bars inserted immediately after $\gamma $ up
until the topological descents all must be inserted left to right to
avoid descents.  The bar insertion which is shifted into $\gamma $ in the
gradient path from $\tau $ to $\sigma $ must be farther to the left than all
these intermediate bar insertions so that there are rank inversions in 
$\tau $ between the bar insertion being shifted and the bar insertions in
the truncated nontrivial intervals.  This means the bar insertion of 
$\gamma $ shifts across these truncated intervals just as it shifts across
descents.
\EOP

Now consider the critical cells $\sigma $ and $\tau $ of example 
~\ref{critpair}.  The facets $F$ and $G$ agreed up through rank 4 which
means that only ranks of 5 or greater may be deleted in a gradient path
from $\tau $ to $\sigma $.  The restriction of $F$ to ranks above its last
nontrivial minimal skipped interval only has a single rank
inversion $(11,12)$ while the restriction of $G$ has inversions
$(11,12), (10,11), (9,11), (8,11), (7,11), (6,11)$ so that the extra 
inversions of $G$ are exactly the inversions of the permutation 
$s_{10}\circ s_9\circ s_8 \circ s_7 \circ s_6 \circ s_5 = 
(6, 7, 8, 9, 10, 11)$.  Any gradient path $d_{i_1}\circ\cdots\circ d_{i_k}$
from $\tau $ to $\sigma $ would need to remove from $G$ 
the same inversions that are
present in the permutation $s_{i_1}\circ\cdots\circ s_{i_k}$ to be a 
legitimate gradient path, but since the cyclic shift $(6, 7, 8, 9, 10, 11)$
has a unique reduced expression there cannot be any other gradient paths from
$\tau $ to $\sigma $.

\subsection{An injective partnership map on lower critical cells}

Section ~\ref{partner} proved that for every lower critical cell there 
exists at least one critical cell with which it could be cancelled.
Now we will choose how to pair critical cells to be cancelled so that 
(1) no two critical cells are paired with
the same critical cell and (2) the pairing is compatible with lexicographic
order (i.e. there is never a critical cell $\gamma $ such that 
$\sigma <_{lex } \gamma <_{lex} \tau $ for $\sigma ,\tau $ a pair of critical
cells we intend to cancel with each other.  Thus each 
lower critical cell is paired with another critical cell such that there is a 
unique gradient path from one to the other both in the lexicographic discrete
Morse function and also in the subsequent discrete Morse functions derived from
it by reversing the unique gradient paths between lexicographically earlier
pairs of matched critical cells.  Our pairing of critical cells may involve
some of the top-dimensional critical cells, but certainly does not need to
involve all of them since our goal is only to cancel all the lower critical
cells.

Recall that for nontrivial minimal 
skipped intervals of types 1 and 3, we checked that given a critical 
cell $\sigma $, then exactly one of the cells obtained by shifting a bar
into or out of the interval will be critical, so we do not need to be 
concerned with the possibility that two different critical cells would 
both be paired with $\sigma $ in this fashion.
In either case, the point was that a bar may be inserted in the last
nontrivial skipped interval of $\tau $ (or some subsequent descent)
to yield a facet with a critical cell if and only the last bar insertion
in the last nontrivial minimal skipped interval (or in the later 
topological descent) may not be extracted.  The if and only if 
statement is determined by whether the block created later immediately 
to the right of the bar to be extracted is bigger or smaller than the 
block immediately to the left of this bar (created when this bar is 
inserted).  Thus we can be sure there is exactly one critical cell that
may be obtained from
$\sigma $ based on the last nontrivial minimal skipped interval in $\sigma $.

Now let us show that whenever two critical cells $\sigma $ and $\tau $ 
are paired as in Section ~\ref{partner}, but there is another
critical cell $\sigma '$ such that $\sigma <_{lex} \sigma' 
<_{lex} \tau $ then $\sigma '$ may be paired with either $\sigma $ or
$\tau $; furthermore, when $\sigma\sim\sigma' $ (resp. $\tau\sim\sigma' $)
then there will also 
exist a critical cell $\rho $
which may be paired with $\tau $ (resp. with $\sigma $).
The critical 
cell $\rho $ exists by virtue of the structure of $\sigma $
after its last nontrivial skipped interval.
We always pair each critical cell with the partner which agrees with
it up to the highest rank; the point is that if there 
is some cell agreeing with $\sigma $ up to a higher rank than $\tau $
does, then the similarity between $\sigma $ and $\tau $ after this
rank ensures
there is also a similar partner critical cell for $\tau $.  

\begin{examp}
\rm
Figure ~\ref{re-pair} gives a saturated chain which has two nontrivial
minimal skipped intervals, each skipping only two ranks.  These 
intervals go from rank 4 to 7 and from rank 14 to 17.  The critical 
cell $\sigma $ of this facet is paired with a lexicographically later
critical cell $\tau $ which has the second nontrivial skipped interval 
of $e\sigma $ replaced by a skipped interval of height one, which means 
Section ~\ref{partner} provides a possible partner
for $\tau $ by exhibiting a different critical cell $\rho $ 
which could be paired with $\tau $ (based on the only nontrivial minimal 
skipped interval of $\tau $,
i.e. the first of two nontrivial minimal skipped intervals for $\sigma $).
One may obtain a partner critical cell for $\rho $, as needed, by 
replacing a topological descent in $e\rho $ with a nontrivial minimal 
skipped interval of height 2 just as is done in converting from $\sigma $
to $\tau $.
\begin{figure}[h]
\begin{picture}(250,70)(40,0)
\includegraphics[width=0.89\textwidth, angle = 0]{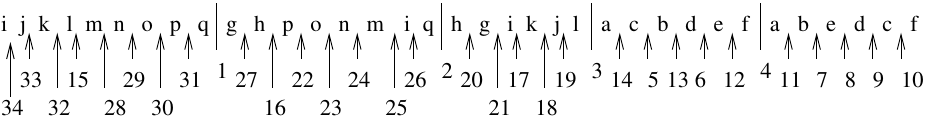}
\end{picture}
\caption{A saturated chain $e\sigma $ to be paired with $\rho $}
\label{re-pair}
\end{figure}
\end{examp}

Suppose the
nontrivial minimal skipped interval of highest rank in $e\sigma'$ 
is an interval skipping two ranks.  Then $\sigma'$ is paired with a critical 
cell $\sigma $ which instead skips either one or three ranks.
In the case of skipping one rank, this is no longer nontrivial, so $\sigma $
will be paired in Section ~\ref{partner} based on an earlier nontrivial
minimal skipped interval with a different critical cell $\rho $.  
However, $\sigma '$ agrees with $\sigma $ to a 
later rank, so we choose $\sigma\sim \sigma '$.  Because $\sigma $ is related
to $\rho $ by the shifting of a single bar insertion to after the nontrivial
minimal skipped interval, $\rho $ possesses the same later structure that 
gave rise to a partner critical cell for $\sigma $, and this structure will
likewise give rise to a partner $\rho '$ for $\rho $.

\subsection{M\"obius functions for $\Pi_n/S_{\lambda }$}

Now that we have proven the  homotopy type for $\lambda $
of hook shape is a wedge of spheres concentrated in top dimension, 
let us conclude with an example of how this may
be helpful for computing M\"obius functions.  

\begin{prop}
Let $\lambda $ be of hook shape.  Then 
$\mu_{P}(\hat{0},\hat{1}) = 0$ if and only if $\lambda_1 > 
3(l(\lambda )-1)$ where $l(\lambda )$ is the length of $\lambda $, i.e.
the number of nonempty blocks.
\end{prop}

\proof
We claim that our lexicographic discrete Morse function for 
$\Pi_n/S_{\lambda }$ has 
uncancelled top-dimensional
critical cells if and only if $\lambda_1 \le 3[l(\lambda )-1]$.
Let $\lambda = (m,1^k)$.  Let us represent elements of $\Pi_n/S_{\lambda }$
by partitions of $\{ a^m,b_1,b_2,\dots ,b_k \} $.
When $m=3,k=1$, then there is a  
critical cell in the facet $a|_3a|_1a|_2b_i $ that 
is not cancelled by any lower critical cell.  More generally, when $m=3k$,
then there is a top-dimensional critical cell that is not cancelled in
the facet which begins by inserting bars left to 
right to form the partition $$aaab_k|aaab_{k-1}|\cdots
|aaab_1$$ and then refines blocks right to left as in the $a^3b_i$ 
example; this critical cell also is not cancelled by any lower one.  
Similar examples may be constructed for $m<3k$.  On the other hand, 
there are no critical cells when $m=4,k=1$ (and similarly for $k=1,m>4$).  
Any saturated chain with 
$m>3k$ will have a topological ascent (and hence no top-dimensional critical
cells) because each refinement preserves the property $m>3k$ for some
component (by the pigeonhole principle), eventually reducing 
to the case $k=1,m\ge 4$.

Since all lower critical cells are cancelled in our proof regarding
the homotopy type for $\lambda $ of hook shape, the nonexistence of 
top dimensional critical cells for $m > 3k$
implies the M\"obius function $\mu_P(\hat{0},
\hat{1})$ must be 0 for such $m,k$. 
\EOP

\subsection{Generalization to arbitrary multisets}

The reason our proof concerning 
the homotopy type of $\Delta (\Pi_n/S_{\lambda})$
does not as it stands apply to all $\lambda $ is that we do not 
know whether every lower critical cell is paired with a cell which is also
critical.  More specifically, the partner cells for critical cells
based on nontrivial minimal skipped intervals of type 2 may or may not 
themselves be critical.  This is because they each have nontrivial
skipped intervals of height one larger than the type 2 interval they are 
replacing, but we do not know 
whether these new intervals are always minimal
skipped intervals.  Thus, we do not know for sure that the saturated 
chains contributing these cells are covered by minimal skipped intervals.
If this were verified, then our proof would apply to all $\lambda $.
  
One way to prove that these intervals must be minimal would be to show that 
a certain puzzle on multisets has no solutions.  One approach might be
to specialize to integer partitions and prove that the resulting puzzle
already has no solutions.
If (conversely) there were a solution in integer
partitions, one could probably translate this to a solution 
in multisets with enough
distinct letters that one could find a lower critical 
cell that cannot be cancelled by any other critical cell, so this would 
probably yield a poset that is not homotopy equivalent to a wedge of spheres
of top dimension.  

The puzzle is as follows: is there a set of integers $n_1,\dots ,n_k$ for
$k\ge 3$ and a more refined set $\{ b_{i,j} | 1\le i \le k, 1\le j \le a_i \}$
of integers so that: 
\begin{enumerate}
\item
$n_i = b_{i,1} + \cdots + b_{i,a_i}$ for $i<k$ with each $a_i\ge 2$.
\item
$n_k = b_{k,1} + b_{k,2} = c_{k,1} + c_{k,2}$ for some $b_{k,1} < b_{k,2}$, 
$c_{k,1} < c_{k,2}$ with $b_{k,1} > c_{k,1}$
\item
The collection of numbers $\{ b_{r,s} | s<k \} \cup \{ c_{k,1} , c_{k,2} \}$  
may be partitioned into blocks $B_1,\dots , B_k$ so that the sum of the 
numbers in block $B_i $ is $n_i$, and the smallest block assigned to
$B_i$, denoted $b_{i,1}'$, is at least as small as $b_{i,1}$.  If 
$b_{i,1}'= b_{i,1}$, then list the blocks $b_{i,j}'$ assigned 
to $B_i$ in increasing order, and we require that the 
first $b_{i,j}'$ to differ from the corresponding block $b_{i,j}$ be 
smaller than $b_{i,j}$.
\end{enumerate}

We are not sure whether the integer partitions puzzle has solutions, 
and suggest it as a possible
way of getting at the homotopy type of $\Delta(\Pi_n /S_{\lambda })$ for
general $\lambda $.  It might even suffice to show that there are no
solutions when the following additional requirement is also imposed 
on the puzzle: that the numbers in 
the set $\{ b_{r,s} | s\le k \} $ which are larger than one
are all distinct.  Likewise, require the 
numbers in the set $\{ b_{r,s} | s < k\}
\cup \{ c_{k,1} ,c_{k,2} \} $ which are larger than one all to be distinct.

\section*{Acknowledgments}

We thank Xun Dong, Art Duval and Robin Forman for proofreading
portions of our paper and providing valuable feedback.

\end{document}